\documentclass[11pt,leqno]{article}

\usepackage{graphicx}
\usepackage{latexsym} 
\usepackage{amsmath}
\usepackage{amssymb}
\usepackage{amsfonts}
\usepackage{enumerate}

\usepackage{esint}
\usepackage{MnSymbol}
\usepackage{theorem}

\theoremstyle{plain}
\newtheorem{teo}{Theorem}[section]

\newtheorem{corollary}[teo]{Corollary}

{\theorembodyfont{\rmfamily}

\newtheorem{remark}[teo]{Remark}

}

\renewcommand{\eqref}[1]{\textnormal{(\ref{#1})}}

\numberwithin{equation}{section}

\newcommand{\cvd}{\hfill$\square$}
\newcommand{\proof}[1]{\noindent\textsc{Proof#1}}

\newcommand{\rmd}{\mathrm{d}}

\title{A Friedrichs-Maz'ya inequality for functions of bounded variation}

\author{Luca Rondi\thanks{Dipartimento di Matematica e Geoscienze,
Universit\`a degli Studi di Trieste, via Valerio, 12/1
34127 Trieste
ITALY. E-mail: \texttt{rondi@units.it}} }

\date{}

\begin{document}

\maketitle

\setcounter{section}{0}
\setcounter{secnumdepth}{2}

\begin{abstract}
The aim of this short note is to give an alternative proof, which applies to functions of bounded variation in arbitrary domains, of an inequality by Maz'ya that improves Friedrichs inequality. A remarkable feature of such a proof is that it is rather elementary, if the basic background in the theory of functions of bounded variation is assumed. Never the less, it allows to extend all the previously known versions of this fundamental inequality to a completely general version. In fact the inequality presented here is optimal in several respects.

As already observed in previous proofs, the crucial step is to provide conditions under which a function of bounded variation on a bounded open set, extended to zero outside, has bounded variation on the whole space. We push such conditions to their limits. In fact, we give a sufficient and necessary condition if the open set has a boundary with $\sigma$-finite surface measure and a sufficient condition if the open set is fully arbitrary.
Via a counterexample we show that such a general sufficient condition is sharp.

\medskip


\medskip

\end{abstract}

\section{Introduction}\label{intro}%

The aim of this paper is to prove by almost elementary methods an extremely general version of a Maz'ya inequality for functions of bounded variation defined on a bounded open set $\Omega\subset \mathbb{R}^N$, $N\geq 2$. Such an inequality implies and extends the so-called Friedrichs inequality. 

We shall use only the basic theory of functions of bounded variation, which may be found for instance in the books \cite{A-F-P,G}, possibly with the single exception of a deep result by Federer, Theorem~4.5.11 in \cite{F}.

Given $\Omega\subset \mathbb{R}^N$, $N\geq 2$, open and bounded, and $u$ a function defined on $\overline{\Omega}$, Friedrichs proved, under suitable assumptions on $\Omega$ and $u$, the following inequality \cite{Fri}
\begin{equation}\label{Friinitial}
\|u\|_{L^2(\Omega)}\leq C\left[\|\nabla u\|_{L^2(\Omega)}+\|u\|_{L^2(\partial\Omega)}\right].
\end{equation}
where the constant $C$ is clearly independent on $u$.

Such an estimate has been greatly improved by Maz'ya in \cite{M0}. He proved
that,
given $\Omega\subset \mathbb{R}^N$, $N\geq 2$, open and bounded, for any $u\in C(\overline{\Omega})\cap W^{1,1}(\Omega)$ we have
\begin{equation}\label{fri0}
\|u\|_{L^{N/(N-1)}(\Omega)}\leq C(N)\left[\|\nabla u\|_{L^1(\Omega)}+\|u\|_{L^1(\partial\Omega)}\right].
\end{equation}
Here and in \eqref{Friinitial}, the last integral is with respect to the surface measure $\mathcal{H}^{N-1}$ and in the sequel of the paper any integral or $L^p$ space over $\partial \Omega$ will be intended 
with respect to the surface measure $\mathcal{H}^{N-1}$, unless explicitly stated otherwise. The constant $C(N)$ in \eqref{fri0}
depends on $N$ only and coincides with the one of the isoperimetric inequality, that is
\begin{equation}\label{constant}
C(N)=\frac{|B_1|^{(N-1)/N}}{\mathcal{H}^{N-1}(\partial B_1)}
\end{equation}
where $B_1$ is any ball of radius $1$ in $\mathbb{R}^N$, 
$|B_1|$ is its $N$-dimensional Lebesgue measure and $\mathcal{H}^{N-1}(\partial B_1)$ is the surface measure of the unit sphere in $\mathbb{R}^N$. In fact the constant $C(N)$ may not be improved.
A proof of this inequality may be found in the book of the same author \cite{M}, see the Corollary on page~319. We shall refer to the inequality \eqref{fri0} as Maz'ya inequality
and we notice that it implies the classical Friedrichs inequality \eqref{Friinitial}, see Remark~\ref{frirem}.

Concerning the corresponding Maz'ya inequality for $BV$ functions, this may be based on the following kind of extension result. Given $u$, a function of bounded variation on a bounded open set $\Omega$, which is extended to zero outside $\Omega$, we look for a condition that allows us to say that $u$ is still of bounded variation on the whole $\mathbb{R}^N$. 

Under the assumption that
$\Omega$ is a set of finite perimeter and $\mathcal{H}^{N-1}(\partial\Omega\backslash\partial^{\ast}\Omega)=0$, where $\partial^{\ast}\Omega$ is the essential boundary of $\Omega$, the extension problem is solved in \cite{M}, see the Lemma on page~496, Section~9.5.5. A corresponding Maz'ya inequality is given in
the same book \cite{M}, in Section~9.5.7, see Theorem~1 on page~499.
We wish to remark that both these results have been extended in \cite{Bu-K} to the case in which $\Omega$ is a bounded open set with finite perimeter and $\partial\Omega$ is countably $\mathcal{H}^{N-1}$-rectifiable.%, that is $\partial\Omega$ is $\mathcal{H}^{N-1}$ rectifiable as in Definition~2.57 in \cite{A-F-P}.

We generalise these results into two directions. We first prove a necessary and sufficient condition for the extension problem for $BV$ functions, under the only assumption that $\Omega$ is a bounded open set such that 
$\partial\Omega$ is $\sigma$-finite with respect to the $\mathcal{H}^{N-1}$ measure, see Theorem~\ref{easycase}, and obtain the corresponding Maz'ya inequality. This result extends the ones in \cite{Bu-K}.

Second, we drop any assumption on the bounded open set $\Omega$ and we obtain our main result, Theorem~\ref{mainresult}, which contains
a most general version of Maz'ya inequality where the bounded open set $\Omega$ is arbitrary and $u$ is just a function of bounded variation on $\Omega$. Finally, in Remark~\ref{continuouscaserem} we show the sharpness of this general result.

Finally, we wish to point out that this result fits into an interesting ongoing research that aims to extend and generalise classical inequalities to arbitrary domains. For example, recently this has been done for Sobolev inequalities in \cite{C-M}.

%$\Omega$ set of finite perimeter and $\mathcal{H}^{N-1}(\partial\Omega\backslash\partial^{\ast}\Omega)=0$.
% Nota: per Maz'ya $\partial^{\ast}E=J_{\chi_E}\cap\partial E$, quindi
%\partial^{\ast}\Omega=J_{\chi_\Omega}\subset\partial\Omega$, vedere pagina 467, Definition 1 e 2. Per gli insiemi aperti $\Omega$ di perimetro finito la condizione
%$\mathcal{H}^{N-1}(\partial\Omega\backslash\partial^{\ast}\Omega)=0$ non dipende dalla definizione di $\partial^{\ast}\Omega$. Per Maz'ya \`e equivalente a dire che esiste, tranne un insieme trascurabile, la normale.

\section{The general Maz'ya inequality}

Throughout the paper, the integer $N\geq 2$ will denote the space dimension. By equivalence with respect to perturbations on sets of measure zero, any measurable set contained in $\mathbb{R}^N$ will be assumed to be a Borel set and any measurable function $u$ which is finite almost everywhere (with respect to the Lebesgue measure) will be assumed to be a real-valued Borel function. For any $s\geq 0$ we denote with $\mathcal{H}^s$ the $s$-dimensional Hausdorff measure.

Let $E\subset \mathbb{R}^N$ be a Borel set. For any $x\in \mathbb{R}^N$ and any $t\in[0,1]$, we say that $E$ has density $t$ at $x$ if
\[\lim_{\rho\to0^+}\fint_{B_{\rho}(x)} \chi_E(y)\rmd y=t\]
where $\chi_E$ denotes the characteristic function of $E$. For any $t\in[0,1]$ we call $E_t$ the subset of points $x$ of $\mathbb{R}^N$ such that $E$ has density $t$ at $x$. We notice that $E_t$ is a Borel set for any $t\in [0,1]$ and that $E_0$ and $E_1$ may be considered as the measure theoretic exterior and interior of $E$, respectively. We call $\partial^{\ast} E$, the essential boundary of $E$, the set 
$\mathbb{R}^N\backslash(E_0\cup E_1)$.

Throughout the paper we fix $\Omega\subset\mathbb{R}^N$, $N\geq 2$, open, and a Borel function $u:\Omega\to \mathbb{R}$.
We shall always tacitly assume that $u$ is extended
to zero outside $\Omega$, that is $u:\mathbb{R}^N\to\mathbb{R}$ is a Borel function  
which is $0$ almost everywhere outside $\Omega$.
Further assumptions on $\Omega$ and $u$ will be specified as needed.

We call $u_+=u\vee 0=\max\{u,0\}$ and $u_-=- (u\wedge 0)=-(\min\{u,0\})$
so that $u=u_+-u_-$ and $|u|=u_++u_-$. Furthermore, for any $t>0$ we define $u_t$ as the following truncated function
$u_t=(u\wedge t)\vee (-t)$.

There are several ways to define the discontinuity set of a real-valued Borel function. We follow Definition 1.57 in \cite{B} and Definitions 3.63 and 3.67 in \cite{A-F-P} and we state a few of their basic properties, see the cited references for details.

For any $x\in \mathbb{R}^N$ we define the approximate upper and lower limits of $u$ at $x$ as follows
\[\begin{array}{l}
\displaystyle{u^+(x)=\mathrm{ap\text{-}}\!\limsup_{y\to x}u(y)=\inf\{t:\ \{u>t\}\text{ has density }0\text{ in }x\}};\\
\displaystyle{u^-(x)=\mathrm{ap\text{-}}\!\liminf_{y\to x}u(y)=\sup\{t:\ \{u>t\}\text{ has density }1\text{ in }x\}}.
\end{array}\]
Clearly $u^-$ and $u^+$ are Borel functions with values in the extended real line and 
$-\infty\leq u^-(x)\leq u^+(x)\leq +\infty$. 
If $u^+(x)=u^-(x)$ we have $u^+(x)-u^-(x)=0$ and we call the common value $\tilde{u}(x)$, the 
approximate limit of $u$ at $x$, that is 
\[\tilde{u}(x)=\displaystyle{u^+(x)=u^-(x)=\mathrm{ap\text{-}}\!\lim_{y\to x}}u(y),\]
and we say that $u$ is approximately continuous at $x$. We call $S(u)$ the set of $x\in\mathbb{R}^N$ such that $u$ is not approximately continuous at $x$.

For any $\Omega\subset\mathbb{R}^N$ open and any function $u\in L^1_{loc}(\Omega)$, we say that $x\in\Omega$ is a Lebesgue point for $u$ if there exists $\tilde{u}(x)\in \mathbb{R}$ such that
\[\lim_{\rho\to 0^+}\fint_{B_{\rho}(x)}|u(y)-\tilde{u}(x)|\rmd y=0.\]
We call $S_u$ the set of points $x\in\Omega$ such that $x$ is not a Lebesgue point for $u$. Since for any Lebesgue point $x$ we have that $u$ is approximately continuous at $x$ and its approximate limit coincide with $\tilde{u}(x)$ we deduce that $S(u)\subset S_u$. Moreover, if $x$ is a Lebesgue point for $u$, it is a Lebesgue point for $|u|$ too and $\widetilde{|u|}(x)=|\tilde{u}|(x)$, thus $S_{|u|}\subset S_u$.
We have that $S_u$ is a Borel set of measure zero
and $\tilde{u}$ is a real-valued Borel function defined on $\Omega\backslash S_u$
that coincides almost everywhere with $u$.
Finally, if $u\in L^{\infty}_{loc}(\Omega)$, then actually $S(u)=S_u$.

Furthermore, for any real-valued Borel function $u$, and any $0<s<t$, we have that
$S(u_s)=S_{u_s}\subset S(u_t)=S_{u_t}\subset S(u)$ and for any $x\in S(u)$ we have 
\[|u_s^+(x)-u_s^-(x)|\leq |u_t^+(x)-u_t^-(x)|\leq |u^+(x)-u^-(x)|.\]
We also have that $S(|u|)$, $S(u_+)$ and $S(u_-)$ are contained in $S(u)$ and for any $x\in S(u)$ we have 
\begin{equation}\label{important}
\left||u|^+(x)-|u|^-(x)\right|
\leq |u^+(x)-u^-(x)|\quad\text{and}\quad\left||u|^+(x)-|u|^-(x)\right| \leq |u|^+(x)
\end{equation}
and
\begin{equation}\label{important1}
\max\left\{\left|u_+^+(x)-u_+^-(x)\right|, \left|u_-^+(x)-u_-^-(x)\right|\right\}
\leq |u^+(x)-u^-(x)|,
\end{equation}
where we used the notation $u_+^+=(u_+)^+$.
Finally,
\[S(u)=\bigcup_{n}S(u_n)\]
from which we deduce that $S(u)$ is a Borel set of measure zero. 

Finally, for any $u\in L^1_{loc}(\Omega)$ and any $x\in \Omega$ we say that $x$ is an approximate jump point for $u$ if there exists a triple $(\tilde{u}^-(x),\tilde{u}^+(x),\nu(x))$,
or equivalently  $(\tilde{u}^+(x),\tilde{u}^-(x),-\nu(x))$, such that 
$\tilde{u}^-(x)$ and $\tilde{u}^+(x)$ are different real numbers, 
 $\nu(x)$ is a unit vector in $\mathbb{R}^N$ and
\[\lim_{\rho\to 0^+}\fint_{B^+_{\rho}(x,\nu(x))}|u(y)-\tilde{u}^+(x)|\rmd y=0\quad\text{and}\quad \lim_{\rho\to 0^+}\fint_{B^-_{\rho}(x,\nu(x))}|u(y)-\tilde{u}^-(x)|\rmd y=0\]
where, for any $\rho>0$, $B^{\pm}_{\rho}(x,\nu(x))=\{y\in B_{\rho}(x): y\cdot\nu(x)\gtrless 0\}$. We denote with $J_u$ the set of approximate jump points of $u$ or jump set of $u$. 
In this case, neither $J_{|u|}\subset J_u$ nor $J_u \subset J_{|u|}$ hold, but for any $x\in J_u$ the previous limits hold for $|u|$ with
$\widetilde{|u|}^+(x)=|\tilde{u}^+(x)|$ and $\widetilde{|u|}^-(x)=|\tilde{u}^-(x)|$.
We notice that $J_u\subset S(u)$ and for any $x\in J_u$ we have
$\tilde{u}^-(x)=u^-(x)$ and $\tilde{u}^+(x)=u^+(x)$, provided we choose $\nu(x)$ such that
$-\infty<\tilde{u}^-(x)<\tilde{u}^+(x)<+\infty$.
 
For any $\Omega\subset\mathbb{R}^N$ open and any function $u\in L^1_{loc}(\Omega)$ we call
the total variation of $u$ on $\Omega$, $|Du|(\Omega)$, the following 
\[|Du|(\Omega)=\sup\left\{\int_{\Omega}u\ \mathrm{div} g:\ g\in C^1_0(\Omega,\mathbb{R}^N)\text{ such that }\|g\|_{L^{\infty}(\Omega)}\leq 1\right\}.\] 
If $|Du|(\Omega)$ is finite we say that $u$ is a function of bounded variation in $\Omega$. We call $BV(\Omega)=\{u\in L^1(\Omega):\ |Du|(\Omega)<+\infty\}$.
%Finally, if $u:\Omega\to \mathbb{R}$ is a Borel function such that $u\mathrm{div}$ for some $g\in C^1_0(\Omega,\mathbb{R}^N)$ is not integrable in $\Omega$, then we also say that $|Du|(\Omega)=+\infty$.

It is a well-known fact that if $u\in L^1_{loc}(\Omega)$ has bounded variation then also 
$u_+$, $u_-$ and $|u|$ have.
Moreover,
\begin{equation}\label{equaz1}
|Du_+|(\Omega)+|Du_-|(\Omega)=|Du|(\Omega)\geq \left|D|u|\right|(\Omega).
\end{equation}
Therefore, if $|Du|(\Omega)$ is finite, for any $t>0$ we have that $u_t\in BV(\Omega)\cap L^{\infty}(\Omega)$ and
 $|Du_t|(\Omega)\leq |Du|(\Omega)$.

We shall also use the following formula.
For any $t\in\mathbb{R}$ we define $\Omega^t=\{x\in\Omega:\ u(x)>t\}\subset\Omega$.
We recall that 
\begin{equation}\label{fleming}
|Du|(\Omega)=\int_{-\infty}^{+\infty}P(\Omega^t,\Omega)\rmd t,
\end{equation}
where $P(\Omega^t,\Omega)=|D\chi_{\Omega^t}|(\Omega)$.

The basic inequality is the following well-known Sobolev inequality for $BV$ functions, see for instance \cite[Theorem~1.28]{G}

\begin{teo}\label{basicestimate}
Let $u\in BV(\mathbb{R}^N)$ be such that $u=0$ almost everywhere outside a bounded set. Then
\[\|u\|_{L^{N/(N-1)}(\mathbb{R}^N)}\leq C(N)|Du|(\mathbb{R}^N).\]
\end{teo}

Notice that $C(N)$ is a constant depending on $N$ only  and it is given by \eqref{constant}. It coincides with the best constant in the isoperimetric inequality, which by the way is an easy consequence of this result, thus it may not be improved.

The inequality we shall prove is the following.
\begin{teo}\label{mainresult}
Let $\Omega\subset\mathbb{R}^N$ be a bounded open set and $u\in L^1_{loc}(\Omega)$, $u$ extended to zero outside $\Omega$. Then
\begin{equation}\label{fri}
\|u\|_{L^{N/(N-1)}(\Omega)}\leq C(N)\left[\left|D|u|\right|(\Omega)+\int_{\partial\Omega}|u|^+(y)\rmd \mathcal{H}^{N-1}(y)\right].
\end{equation}
\end{teo}

\begin{remark}\label{omega0rem}
Let us notice that for any $x\in \partial\Omega\cap\Omega_0$, then, no matter what Borel function $u:\Omega\to\mathbb{R}$ we have, $u$ and $|u|$ are approximately continuous at $x$ and $\tilde{u}(x)=\widetilde{|u|}(x)=0$, therefore $|u|^+(x)=0$ as well. Moreover, $S(u)\cap \partial\Omega\cap\Omega_0=S(|u|)\cap \partial\Omega\cap\Omega_0=\emptyset.$
\end{remark}

We notice that we have essentially no assumption on $u$ and $\Omega$. Of course the result is trivial if the right hand side is equal to $+\infty$.

Let us mention and prove the following corollary, which is a slightly more general version of the result of Corollary 1, page~391, in \cite{M}. 

\begin{corollary}\label{rgeq1cor}
Let $\Omega\subset\mathbb{R}^N$ be a bounded open set and $u\in L^1_{loc}(\Omega)$, $u$ extended to zero outside $\Omega$.

Let us fix
\begin{equation}\label{rpq}
r\geq 1,\quad q=\frac{rN}{N-1}\geq \frac{N}{N-1},\quad 1\leq p=\frac{rN}{N-1+r} <N
\end{equation}
or, equivalently,
\[1\leq p<N,\quad r=\frac{p(N-1)}{N-p}\geq 1,\quad q=\frac{pN}{N-p}.\]

Assume that $\nabla u\in L^p(\Omega,\mathbb{R}^N)$, then there exists a constant $C(N,r)$, depending on $N$ and $r$ only, such that
\begin{equation}\label{frir}
\|u\|_{L^q(\Omega)}\leq C(N,r)\left[\|\nabla u\|_{L^p(\Omega)}+\||u|^+\|_{L^r(\partial\Omega)}\right].
\end{equation}
\end{corollary}

\begin{remark}\label{frirem} We notice that from \eqref{frir} we can easily deduce the classical Friedrichs inequality \eqref{Friinitial}. In fact, if we choose $r=2$, then
$q=2N/(N-1)>2$ and $p=2N/(N+1)<2$. Therefore a simple application of H\"older inequality leads to \eqref{Friinitial} with a constant $C$ depending on $N$ and $|\Omega|$ only.
\end{remark}

\proof{.} The case $r=1$ is included in Theorem~\ref{mainresult}, thus we assume $r>1$.
For the time being, let us assume that $u\in L^{\infty}(\Omega)$. 
We take $|u|^r$ and we observe that $(|u|^r)^+=(|u|^+)^r$ on $\partial\Omega$ and that
\[\int_{\Omega}\|\nabla (|u|^r)\|=\int_{\Omega}r |u|^{r-1}\|\nabla u\|\leq
r\left(\int_{\Omega}|u|^{rN/(N-1)}\right)^{(N-1)(r-1)/(rN)}\|\nabla u\|_{L^p(\Omega)}.\]
We can apply Theorem~\ref{mainresult} to $|u|^r$ and obtain that
\[\|u\|_{L^q(\Omega)}^r\leq C(N)\left[r\|u\|_{L^q(\Omega)}^{r-1}\|\nabla u\|_{L^p(\Omega)}+
\||u|^+\|_{L^r(\partial\Omega)}^r
\right].\]
For any constant $\lambda>0$, we have that
\[r\|u\|_{L^q(\Omega)}^{r-1}\|\nabla u\|_{L^p(\Omega)}\leq (r-1)\lambda^{r/(r-1)}\|u\|_{L^q(\Omega)}^r+\lambda^{-r}\|\nabla u\|_{L^p(\Omega)}^r.\]
Choosing $\lambda$ such that
$C(N)(r-1)\lambda^{r/(r-1)}=1/2$, we infer that
\[\|u\|_{L^q(\Omega)}^r\leq 2C(N)\left[\lambda^{-r}
\|\nabla u\|_{L^p(\Omega)}^r+
\||u|^+\|_{L^r(\partial\Omega)}^r
\right],\]
that is
\[\|u\|_{L^q(\Omega)}\leq (2C(N))^{1/r}\left[\lambda^{-1}
\|\nabla u\|_{L^p(\Omega)}+
\||u|^+\|_{L^r(\partial\Omega)}
\right].\]
An easy computation shows that
\[(2C(N))^{1/r}\max\{\lambda^{-1},1\}\leq 2C(N)r.\]
and the proof is concluded by choosing $C(N,r)=2C(N)r$.

We can easily drop the assumption that $u\in L^{\infty}(\Omega)$ by taking $u_t$ with $t>0$ and letting $t\to+\infty$.\cvd
%\bigskip

Clearly, the main issue to solve in order to prove Theorem~\ref{mainresult} is the following. Assuming that $\Omega$ is a bounded open set and $u$ is a function of bounded variation in $\Omega$, which conditions are sufficient to have that $u$, extended to $0$ outside $\Omega$, belongs to $BV(\mathbb{R}^N)$? And in this case, what is the relation between $|Du|(\mathbb{R}^N)$ and $|Du|(\Omega)$?

Before proving 
Theorem~\ref{mainresult} we shall state and prove a sharper result under the assumption that the boundary of $\Omega$ is $\sigma$-finite with respect to the $\mathcal{H}^{N-1}$ measure. We recall that, by a $\sigma$-finite set with respect to a measure, we mean a measurable set which may be obtained as the union of a sequence of measurable sets with finite measure.

\begin{teo}\label{easycase}
Let $\Omega\subset \mathbb{R}^N$ be a bounded open set such that 
$\partial\Omega$ is $\sigma$-finite with respect to the $\mathcal{H}^{N-1}$ measure.
Let $u$ be a function of bounded variation on $\Omega$, that is $u\in L^1_{loc}(\Omega)$ such that $|Du|(\Omega)<+\infty$.

Then $u\in BV(\mathbb{R}^N)$ if and only if 
\begin{equation}\label{iff}
\int_{\partial\Omega\cap S(u)}\left|u^+(y)-u^-(y)\right|\rmd\mathcal{H}^{N-1}(y)<+\infty.
\end{equation}
In this case
\begin{equation}\label{firstequation}
|Du|(\mathbb{R}^N)=|Du|(\Omega)+\int_{\partial\Omega\cap S(u)}\left|u^+(y)-u^-(y)\right|\rmd\mathcal{H}^{N-1}(y)
\end{equation}
and
\begin{equation}\label{frieasy}
\|u\|_{L^{N/(N-1)}(\Omega)}\leq C(N)\left[\left|D|u|\right|(\Omega)+
\int_{\partial\Omega\cap S(|u|)}\left||u|^+(y)-|u|^-(y)\right|\rmd\mathcal{H}^{N-1}(y)\right].
\end{equation}
\end{teo}

\begin{remark}\label{easycaserem}
Let us observe that  in all the integrals above we can replace $\partial\Omega\cap S(u)$ and $\partial\Omega\cap S(|u|)$
with $\partial\Omega$ or, by Remark~\ref{omega0rem},  with
$\partial\Omega\cap S(u)\backslash \Omega_0$ and $\partial\Omega\cap S(|u|)\backslash \Omega_0$, respectively.
Moreover we have that, by \eqref{equaz1},
$\left|D|u|\right|(\Omega)\leq |Du|(\Omega)$
and, by \eqref{important},
\begin{multline}\label{equaz2}
\int_{\partial\Omega\cap S(|u|)}\left||u|^+(y)-|u|^-(y)\right|\rmd\mathcal{H}^{N-1}(y)\leq
\int_{\partial\Omega\cap S(u)}\left|u^+(y)-u^-(y)\right|\rmd\mathcal{H}^{N-1}(y)\quad\text{and}\\
\int_{\partial\Omega\cap S(|u|)}\left||u|^+(y)-|u|^-(y)\right|\rmd\mathcal{H}^{N-1}(y)\leq\int_{\partial\Omega}|u|^+(y)\rmd\mathcal{H}^{N-1}(y).
\end{multline}
Thus we have obtained, under this assumption on $\Omega$, a perfectly sharp and improved form of Theorem~\ref{mainresult}. Let us finally notice that we do not even require that $\Omega$ has finite perimeter.
\end{remark}

\begin{remark}\label{easycasecor2}
Also Corollary~\ref{rgeq1cor} may be improved in this case. In fact, assume that $\Omega\subset \mathbb{R}^N$ is a bounded open set such that 
$\partial\Omega$ is $\sigma$-finite with respect to the $\mathcal{H}^{N-1}$ measure,
thus $\Omega$ may not have finite perimeter.
Let $u\in L^1_{loc}(\Omega)$, $u$ extended to zero outside $\Omega$. Let $r$, $p$ and $q$ be as in \eqref{rpq} and assume that $\nabla u\in L^p(\Omega,\mathbb{R}^N)$. Then, for the same constant $C(N,r)$ appearing in \eqref{frir}, we have
\begin{equation}\label{frir2}
\|u\|_{L^q(\Omega)}\leq C(N,r)\left[\|\nabla u\|_{L^p(\Omega)}+
\left(\int_{\partial\Omega\cap S(|u|)}\left|(|u|^+)^r(y)-(|u|^-)^r(y)\right|\rmd\mathcal{H}^{N-1}(y)\right)^{1/r}\right].
\end{equation}
\end{remark}

\proof{.}
One implication is easy. If $u\in BV(\mathbb{R}^N)$ then
$|Du|(\mathbb{R}^N)=|Du|(\Omega)+|Du|(\partial\Omega)+|Du|(\mathbb{R}^N\backslash\overline{\Omega})=|Du|(\Omega)+|Du|(\partial\Omega)$.

By Theorems~3.78 and 3.77 and Lemma 3.76 in \cite{A-F-P}, $\mathcal{H}^{N-1}(S_u\backslash J_u)=0$ and
\begin{multline}\label{firstformula}
|Du|(\partial\Omega)=|Du|(\partial\Omega\cap S_u)=\int_{\partial\Omega\cap S_u}\left|u^+(y)-u^-(y)\right|\rmd\mathcal{H}^{N-1}(y)=\\\int_{\partial\Omega\cap S(u)}\left|u^+(y)-u^-(y)\right|\rmd\mathcal{H}^{N-1}(y)=
\int_{\partial\Omega\cap J_u}\left|u^+(y)-u^-(y)\right|\rmd\mathcal{H}^{N-1}(y).
\end{multline}
Here the assumption that $\partial\Omega$ is $\sigma$-finite with respect to the $\mathcal{H}^{N-1}$ measure is used to be sure that, on $\partial\Omega$,
$Du$ coincides with its jump part $D^j(u)$ and no contribution is due by the
absolutely continuous part and by the Cantor part of $Du$.
Since $|Du|(\partial\Omega)\leq |Du|(\mathbb{R}^N)<+\infty$, we have that \eqref{iff} holds. Moreover, \eqref{firstequation} holds as well, whereas \eqref{frieasy}
immediately follows from \eqref{firstequation} and Theorem~\ref{basicestimate}.

%$Du\lefthalfcup $

We now deal with the other implication. We divide the proof into three cases.

\smallskip

\noindent
\textsc{First case.} Assume that $\mathcal{H}^{N-1}(\partial\Omega)<+\infty$ and that $u\in L^{\infty}(\Omega)$. We claim that if $u$ has bounded variation on $\Omega$, then $u\in BV(\mathbb{R}^N)$ and, by the previous part of the proof, \eqref{firstequation} holds.

The proof of this claim follows the proof of Proposition 3.62 in \cite{A-F-P} that shows that, under these assumptions, $\Omega$ is a set of finite perimeter. There exists a constant $C_1(N)$, depending on $N$ only, such that for any $\delta$, $0<\delta<1$, we can find $x_i\in\partial\Omega$ and $r_i>0$, $i=1,\ldots,n$, such that for any $i=1,\ldots,n$ we have
$2r_i<\delta$ and, setting $B_i^{\delta}=B_{2r_i}(x_i)$,
\[\partial\Omega\subset B^{\delta}=\bigcup_{i=1}^nB_i^{\delta}\quad\text{and}\quad
\sum_{i=1}^n\mathcal{H}^{N-1}(\partial B_i^{\delta})\leq C_1(N)\left(\mathcal{H}^{N-1}(\partial\Omega)+\delta\right).\]

We call $u^{\delta}$ the function which is equal to $u$ in $\Omega\backslash \overline{B^{\delta}}$ and zero otherwise. An easy application of Theorem~3.84 in \cite{A-F-P} shows that $u^{\delta}\in BV(\mathbb{R}^N)$. Moreover, for any $\delta$, $0<\delta<1$,
\begin{multline*}
|Du^{\delta}|(\mathbb{R}^N)\leq |Du|(\Omega\backslash \overline{B^{\delta}})+
2\|u\|_{L^{\infty}(\mathbb{R}^N)}\sum_{i=1}^n\mathcal{H}^{N-1}(\partial B_i^{\delta})\leq\\
|Du|(\Omega)+
2\|u\|_{L^{\infty}(\mathbb{R}^N)}C_1(N)\left(\mathcal{H}^{N-1}(\partial\Omega)+1\right).
\end{multline*}

Since, as $\delta\to 0^+$, $u^{\delta}$ converges to $u$ in $L^1(\mathbb{R}^N)$ we easily conclude that $u\in BV(\mathbb{R}^N)$ as well.

\smallskip

\noindent
\textsc{Second case.} We drop the assumption that $\mathcal{H}^{N-1}(\partial\Omega)<+\infty$. However we assume that $u\in L^{\infty}(\Omega)$ and that $u\geq 0$ everywhere in $\Omega$.%, by replacing $u$ with $|u|$ if needed.

For any $t\in\mathbb{R}$ we define $\Omega^t=\{x\in\Omega:\ u(x)>t\}\subset\Omega$ as before.
We wish to show that 
$|Du|(\mathbb{R}^N)=\int_{-\infty}^{+\infty}P(\Omega^t,\mathbb{R}^N)\rmd t$ is finite.
We remark that since $u\geq 0$ we can restrict this integral and the one in \eqref{fleming} to the interval $(0,+\infty)$.

We have that, for any $t>0$,
\begin{equation}\label{omegatest}
\{x\in\partial\Omega:\ u^+(x)> t\}
\subset(\partial^{\ast}\Omega^t\cup\Omega^t_1)\cap\partial\Omega\subset \{x\in\partial\Omega:\ u^+(x)\geq t\}.
\end{equation}
%Moreover, for any $0<s<t$ we have
%$(\partial^{\ast}\Omega^t\cup\Omega^t_1)\cap\partial\Omega\subset
%(\partial^{\ast}\Omega^s\cup\Omega^s_1)\cap\partial\Omega$.
Moreover, for any $x\in \partial\Omega$, if $x\in \partial^{\ast}\Omega^t$ then
$u^-(x)\leq t\leq u^+(x)$. For any $x\in\partial\Omega$ we call $E_x=\{t>0:\ x\in \partial^{\ast}\Omega^t\}$ and $|E_x|=\int_{0}^{+\infty}\chi_{E_x}(t)\rmd t$. We notice that
$E_x\subset [u^-(x), u^+(x)]$, hence $|E_x|\leq |u^+(x)-u^-(x)|$.
Therefore, since $\partial\Omega$ is $\sigma$-finite with respect to the $\mathcal{H}^{N-1}$ measure, we can use the Fubini Theorem and obtain
\begin{multline*}
\int_{0}^{+\infty}\mathcal{H}^{N-1}(\partial^{\ast}\Omega^t\cap\partial\Omega)\rmd t
=\int_{\partial\Omega}|E_y|
\rmd\mathcal{H}^{N-1}(y)
\leq\\
\int_{\partial\Omega}\left|u^+(y)-u^-(y)\right|
\rmd\mathcal{H}^{N-1}(y)=
\int_{\partial\Omega\cap S(u)}\left|u^+(y)-u^-(y)\right|\rmd\mathcal{H}^{N-1}(y).
\end{multline*}

We conclude that
\begin{multline}\label{crucial}
\int_{0}^{+\infty}\mathcal{H}^{N-1}(\partial^{\ast}\Omega^t)\rmd t
=\int_{0}^{+\infty}\mathcal{H}^{N-1}(\partial^{\ast}\Omega^t\cap\Omega)\rmd t+
\int_{0}^{+\infty}\mathcal{H}^{N-1}(\partial^{\ast}\Omega^t\cap\partial\Omega)\rmd t\leq\\
|Du|(\Omega)+\int_{\partial\Omega\cap S(u)}\left|u^+(y)-u^-(y)\right|\rmd\mathcal{H}^{N-1}(y).
\end{multline}

Now we use a deep result by Federer, see Theorem~4.5.11 in \cite{F}, that guarantees that for any bounded  Borel set $E$ such that $\mathcal{H}^{N-1}(\partial^{\ast}E)$ is finite, then $E$ is a set of finite perimeter and, consequently, $P(E,\mathbb{R}^N)=\mathcal{H}^{N-1}(\partial^{\ast}E)$.

Hence, for almost any $t>0$, $\mathcal{H}^{N-1}(\partial^{\ast}\Omega^t)$
is finite and coincides with $P(\Omega^t,\mathbb{R}^N)$. Then \eqref{crucial} guarantees that $u\in BV(\mathbb{R}^N)$.

\smallskip

\noindent
\textsc{General case.} For the time being we drop the assumption 
that $u\in L^{\infty}(\Omega)$ but we keep the fact that $u\geq 0$
everywhere in $\Omega$.

However, by the previous steps, for any $t>0$, we have that
$u_t\in BV(\mathbb{R}^N)$, $|Du_t|(\Omega)\leq |Du|(\Omega)$ and
\begin{multline*}|Du_t|(\partial\Omega)=|Du_t|(\partial\Omega\cap S_{u_t})=\int_{\partial\Omega\cap S(u_t)}\left|u_t^+(y)-u_t^-(y)\right|\rmd\mathcal{H}^{N-1}(y)\leq \\
\int_{\partial\Omega\cap S(u_t)}\left|u^+(y)-u^-(y)\right|\rmd\mathcal{H}^{N-1}(y)\leq
\int_{\partial\Omega\cap S(u)}\left|u^+(y)-u^-(y)\right|\rmd\mathcal{H}^{N-1}(y)<+\infty.
\end{multline*}

We conclude that, for some constant $C$ independent of $t>0$, we have 
\[\|u_t\|_{L^{N/(N-1)}(\mathbb{R}^N)}+|Du_t|(\mathbb{R}^N)\leq C\quad\text{for any }t>0.\]
Since, as $t\to+\infty$, $|u_t|$ converges monotonically to $|u|$ everywhere in $\mathbb{R}^N$, we have that
$\|u_t\|_{L^{N/(N-1)}(\mathbb{R}^N)}$ converges to $\|u\|_{L^{N/(N-1)}(\mathbb{R}^N)}$. We conclude that $u\in L^1(\mathbb{R}^N)$ and since, as $t\to+\infty$, $u_t$ converges to $u$ everywhere in $\mathbb{R}^N$ and in $L^1(\mathbb{R}^N)$ we immediately infer that $u\in BV(\mathbb{R}^N)$.

Finally we can easily drop the assumption that $u\geq 0$ everywhere in $\Omega$, by using $u_+$ and $u_-$ and the estimates
\eqref{equaz1} and \eqref{important1}.
The proof is concluded.\cvd

\begin{remark}\label{easycaserem2}
Under the assumptions of Theorem~\ref{easycase}, let us assume that also
 \eqref{iff} holds. Then we can compute $|Du|(\partial\Omega)$ by using the formula in \eqref{firstequation} or the equivalent formulations in \eqref{firstformula}, see also Remark~\ref{easycaserem}.

Another interesting remark is the following. Assume that $\Omega$ is a bounded open set such that $\Omega$ is a set of finite perimeter. % and $\mathcal{H}^{N-1}(\partial\Omega\backslash\partial^{\ast}\Omega)=0$. 
Then \eqref{iff} and \eqref{firstequation} may be replaced, respectively, by
\begin{equation}\label{iff2}
\int_{\partial^{\ast}\Omega}|u|^+(y)\rmd\mathcal{H}^{N-1}(y)+
\int_{\partial\Omega\cap \Omega_1}\left|u^+(y)-u^-(y)\right|\rmd\mathcal{H}^{N-1}(y)
<+\infty
\end{equation}
and, provided \eqref{iff2} holds,
\begin{equation}\label{firstequation2}
|Du|(\mathbb{R}^N)=|Du|(\Omega)+\int_{\partial^{\ast}\Omega}
|u|^+(y)\rmd\mathcal{H}^{N-1}(y)+
\int_{\partial\Omega\cap \Omega_1}\left|u^+(y)-u^-(y)\right|\rmd\mathcal{H}^{N-1}(y).
\end{equation}

%since $J_u\cap\partial\Omega$ and $J_{|u|}\cap\partial\Omega$ differ by an $\mathcal{H}^{N-1}$-negligible set and
%\begin{multline*}
%\int_{\partial\Omega\cap S(u)}
%\left|u^+(y)-u^-(y)\right|\rmd\mathcal{H}^{N-1}(y)=
%\int_{\partial^{\ast}\Omega\cap J_u}\max\{
%|u^+(y)|,|u^-(y)|\}\rmd\mathcal{H}^{N-1}(y)=\\
%\int_{\partial^{\ast}\Omega\cap J_u}
%|u|^+(y)\rmd\mathcal{H}^{N-1}(y)=\int_{\partial\Omega}
%|u|^+(y)\rmd\mathcal{H}^{N-1}(y)=\\
%\int_{\partial\Omega\cap S(|u|)}
%\left||u|^+(y)-|u|^-(y)\right|\rmd\mathcal{H}^{N-1}(y).
%\end{multline*}
In fact in this case, for $\mathcal{H}^{N-1}$-almost any $x\in\partial^{\ast}\Omega$,
$x\in J_{\chi_{\partial\Omega}}$, with triple $(1,0,\nu(x))$, where $\nu(x)$ is, in a  measure theoretic sense, the exterior normal to $\Omega$ at $x$.
Then we can define a Borel function $\tilde{u}^-$ on $\partial^{\ast}\Omega$, with values in the extended real line, such that $\tilde{u}^-(x)=0$ whenever $x\not\in S_u$ and, for 
$\mathcal{H}^{N-1}$-almost any $x\in J_{\chi_{\partial\Omega}}\cap J_u$, $x$ is an approximate jump point for $u$ with triple 
$(\tilde{u}^-(x),0,\nu(x))$, $\nu(x)$ being the exterior normal to $\Omega$ at $x$.
Such a function $\tilde{u}^-$ may be considered as the trace of $u$ on the essential boundary of $\Omega$ and we have
\[|\tilde{u}^-(x)|=\widetilde{|u|}^-(x)=
\left|u^+(x)-u^-(x)\right|=
%\max\{|u^+(x)|,|u^-(x)|\}=
|u|^+(x)\quad\text{ for }\mathcal{H}^{N-1}\text{-almost any }x\in\partial^{\ast}\Omega.\]  
\end{remark}

We now prove the 
Maz'ya inequality in the general case, Theorem~\ref{mainresult}.

\proof{ of Theorem~\ref{mainresult}.}
We conclude the proof of our main result.
Without loss of generality we can assume that
$u\geq 0$ everywhere in $\Omega$, by replacing $u$ with $|u|$ if needed, and that the right hand side of \eqref{fri} is finite.

For the time being we also assume that $u\in L^{\infty}(\Omega)$.
We wish to show that 
$|Du|(\mathbb{R}^N)=\int_{0}^{+\infty}P(\Omega^t,\mathbb{R}^N)\rmd t$ is finite.

We observe that, by the Fubini Theorem, since
$\{y\in\partial\Omega:\ u^+(y)> 0\}$ is $\sigma$-finite with respect to the $\mathcal{H}^{N-1}$ measure,

\[\int_{\partial\Omega}u^+(y)\rmd \mathcal{H}^{N-1}(y)=
\int_0^{+\infty}\mathcal{H}^{N-1}(\{y\in\partial\Omega:\ u^+(y)\geq t\})\rmd t.%=\int_0^{+\infty}\mathcal{H}^{N-1}(\{y\in\partial\Omega:\ u^+(y)>t\})\rmd t.
\]
%We also have that, for any $t>0$,
%\[\{x\in\partial\Omega:\ u^+(x)> t\}
%\subset(\partial^{\ast}\Omega^t\cup\Omega^t_1)\cap\partial\Omega\subset \{x\in\partial\Omega:\ u^+(x)\geq t\}.\] 
%Moreover, for any $0<s<t$ we have
%$(\partial^{\ast}\Omega^t\cup\Omega^t_1)\cap\partial\Omega\subset
%(\partial^{\ast}\Omega^s\cup\Omega^s_1)\cap\partial\Omega$.
By \eqref{omegatest}, we deduce that 
\[t\mathcal{H}^{N-1}((\partial^{\ast}\Omega^t\cup\Omega^t_1)\cap\partial\Omega)
\leq\int_{\partial\Omega}u^+(y)\rmd \mathcal{H}^{N-1}(y).\]
Since
\[\mathcal{H}^{N-1}(\partial^{\ast}\Omega^t)=\mathcal{H}^{N-1}(\partial^{\ast}\Omega^t\cap\Omega)+\mathcal{H}^{N-1}(\partial^{\ast}\Omega^t\cap\partial\Omega)=
P(\Omega^t,\Omega)+\mathcal{H}^{N-1}(\partial^{\ast}\Omega^t\cap\partial\Omega),
\]
we obtain that for almost any $t>0$ we have that $\mathcal{H}^{N-1}(\partial^{\ast}\Omega^t)$ is finite. 
We use again Theorem~4.5.11 in \cite{F} and conclude that, for almost any $t>0$, we have that
\[P(\Omega^t,\mathbb{R}^N)=P(\Omega^t,\Omega)+\mathcal{H}^{N-1}(\partial^{\ast}\Omega^t\cap\partial\Omega)\leq 
P(\Omega^t,\Omega)+\mathcal{H}^{N-1}(\{y\in\partial\Omega:\ u^+(y)\geq t\})\]
therefore
\begin{multline*}
|Du|(\mathbb{R}^N)=\int_0^{+\infty}\!\!P(\Omega^t,\mathbb{R}^N)\rmd t\leq\\
\int_0^{+\infty}\!\!P(\Omega^t,\Omega)\rmd t+
\int_0^{+\infty}\!\!\mathcal{H}^{N-1}(\{y\in\partial\Omega:\ u^+(y)\geq t\})\rmd t
=
|Du|(\Omega)+
\int_{\partial\Omega}u^+(y)\rmd \mathcal{H}^{N-1}(y).
\end{multline*}

If $u$ does not belong to $L^{\infty}(\Omega)$ we argue by the usual truncation argument. The proof is concluded.\cvd
%\bigskip

We conclude the paper with two interesting remarks.
We need a few preliminary considerations. Let $\Omega\subset\mathbb{R}^N$ be open and bounded.
Let us assume that $u:\Omega\to\mathbb{R}$ is the restriction of a continuous function $w:\mathbb{R}^N\to \mathbb{R}$. As usual we extend $u$ to zero outside $\Omega$. 
Obviously $u\in L^{\infty}(\Omega)$ and we have that
$\tilde{u}=u=w$ everywhere in $\Omega$, whereas $\tilde{u}=0$ everywhere in $\mathbb{R}^N\backslash\overline{\Omega}$.

On $\partial\Omega$ we have the following properties. Let us fix $x\in\partial\Omega$.
If $w(x)=0$, then $u$ is approximately continuous at $x$, $|u^+(x)-u^-(x)|=0$ and $\tilde{u}(x)=0$. If $w(x)>0$, then $0\leq u^-(x)\leq u^+(x)\leq w(x)$, whereas if $w(x)<0$, then $w(x)\leq u^-(x)\leq u^+(x)\leq 0$. Overall, we obtain that for any $x\in\partial\Omega$ we have
$|u^+(x)-u^-(x)|\leq |w(x)|$. The value of $u^+(x)$ and $u^-(x)$ may also depend on the density of $\Omega$ at $x$. If $x\in \Omega_0$ then, as already noticed in Remark~\ref{omega0rem}, $u$ is approximately continuous at $x$, $|u^+(x)-u^-(x)|=0$ and $\tilde{u}(x)=0$, no matter what $w(x)$ is. If $x\in\Omega_1$ then 
$u$ is approximately continuous at $x$, $|u^+(x)-u^-(x)|=0$ and $\tilde{u}(x)=w(x)$.
For any $x\in\partial^{\ast}\Omega$ then $u^+(x)=w(x)\vee 0$ and $u^-(x)=w(x)\wedge 0$, therefore $|u^+(x)-u^-(x)|=|w(x)|$.
Therefore, we conclude that 
\begin{multline}\label{continuous}
\int_{\partial\Omega}\left|u^+(y)-u^-(y)\right|\rmd \mathcal{H}^{N-1}(y)=
\int_{\partial^{\ast}\Omega}\left|w(y)\right|\rmd \mathcal{H}^{N-1}(y)=
\int_{\partial^{\ast}\Omega}|u|^+(y)\rmd \mathcal{H}^{N-1}(y)\leq\\
\int_{\partial\Omega}|u|^+(y)\rmd \mathcal{H}^{N-1}(y)=
\int_{\partial\Omega\backslash\Omega_0}|w(y)|\rmd \mathcal{H}^{N-1}(y)
\leq
\int_{\partial\Omega}\left|w(y)\right|\rmd \mathcal{H}^{N-1}(y).
\end{multline}
The first remark is that, consequently, our results clearly include the one in \eqref{fri0}.

In this second remak we prove by an example that Theorem~\ref{mainresult} is sharp. If $\partial\Omega$ is $\sigma$-finite with respect to the $\mathcal{H}^{N-1}$ measure, then Theorem~\ref{easycase} is perfectly sharp. We show that Theorem~\ref{mainresult} is sharp if $\partial\Omega$ is not $\sigma$-finite with respect to the $\mathcal{H}^{N-1}$ measure. 

\begin{remark}\label{continuouscaserem}
Let $C\subset [0,1]$ be the Cantor set and $f:[0,1]\to[0,1]$ be the Cantor function.
We define $g:[0,2]\to[0,1]$ as follows. For any $t\in[0,1]$ we set $g(t)=f(t)$ and for any $t\in[1,2]$ we set $g(t)=f(2-t)$.
We also define $C_1=C\cup \{t=2-s:\ s\in C\}$.

For any $l>0$, let $\Omega^l=\left((0,2)\backslash C_1\right)\times (0,l)\subset\mathbb{R}^2$ and $u^l:\Omega^l\to\mathbb{R}$ be defined as follows
\[u^l(x,y)=g(x)\quad\text{for any }(x,y)\in\Omega^l.\]
Clearly $u^l$ is the restriction to $\Omega^l$ of a continuous function $w:\mathbb{R}^2\to[0,1]$. We also have that $|Du^l|(\Omega^l)=0.$
Furthermore, we have that
$\mathcal{H}^{1}(\partial\Omega^l\backslash(\Omega^l_{1/2}\cup\Omega^l_1))=0$ and, calling $\tilde{\Omega}^l=(0,2)\times (0,l)$ and using \eqref{continuous}, we have
\begin{multline}\label{nonso}
\int_{\partial\Omega^l}\left|(u^l)^+(y)-(u^l)^-(y)\right|\rmd \mathcal{H}^{N-1}(y)=
\int_{\partial\Omega^l\cap S(u^l)}\left|(u^l)^+(y)-(u^l)^-(y)\right|\rmd \mathcal{H}^{1}(y)=\\\int_{\partial^{\ast}\Omega^l}|u^l|^+(y)\rmd \mathcal{H}^{1}(y)=
\int_{\partial\tilde{\Omega}^l}|u^l|^+(y)\rmd \mathcal{H}^{N-1}(y)=2\int_0^2g(t)\rmd t.
\end{multline}

It is easy to show that
$Dw=D^cw$ in $\tilde{\Omega}^l$, where $D^cw$ denotes the Cantor part of $Dw$. Moreover, $Du^l=Dw=0$ in $\Omega^l$. Clearly, the function $u^l$, as usual extended to zero outside $\Omega^l$, coincides almost everywhere with 
$w\chi_{\tilde{\Omega}^l}\in BV(\mathbb{R}^N)$
and, since $|Dw|(\tilde{\Omega}^l)=|D^cw|(\tilde{\Omega}^l)=|D^cw|(\tilde{\Omega}^l\cap \partial\Omega^l)$,
we have
\[|Du^l|(\mathbb{R}^N)=|D(w\chi_{\tilde{\Omega}^l})|(\mathbb{R}^N)=\left(|Dw|(\tilde{\Omega}^l)+2\int_0^2g(t)\rmd t\right)=\left(|D^cw|(\tilde{\Omega}^l\cap \partial\Omega^l)+2\int_0^2g(t)\rmd t\right)\]
and
\[\|u^l\|_{L^{2}(\Omega^l)}=\|w\|_{L^{2}(\tilde{\Omega}^l)}\leq C(2)
\left(|Dw|(\tilde{\Omega}^l)+2\int_0^2g(t)\rmd t\right).\]
We have that $\|w\|_{L^{2}(\tilde{\Omega}^l)}=l^{1/2}\|w\|_{L^{2}(\tilde{\Omega}^1)}$, for any $l>0$, $|Du^l|(\Omega^l)=0$, and $|Dw|(\tilde{\Omega}^l)>0$. Therefore, by \eqref{nonso}, formula
\eqref{firstequation} does not hold for $u^l$ and $\|u^l\|_{L^{2}(\Omega^l)}$ can not be bounded by a constant, independent on $l$, times $|Du^l|(\Omega^l)$ plus either
$\int_{\partial\Omega^l\cap S(u^l)}\left|(u^l)^+(y)-(u^l)^-(y)\right|\rmd \mathcal{H}^{1}(y)$ or 
$\int_{\partial^{\ast}\Omega^l}|u^l|^+(y)\rmd \mathcal{H}^{1}(y)$. It is indeed necessary to add $\int_{\partial\Omega^l\cap\Omega^l_1}|u^l|^+(y)\rmd \mathcal{H}^{N-1}(y)$.
%Let us finally observe that such a last term provides us with a bound for the total variation of $Dw=D^cw$ on $\tilde{\Omega}^l$.
\end{remark}

%\bigskip

\subsubsection*{Acknowledgements}
Luca Rondi is supported by Universit\`a 
degli Studi di Trieste through FRA 2014 and by GNAMPA, INdAM, through 2015 projects.

\end{document}